\documentclass[12pt]{amsart}

\newtheorem{theorem}{Theorem}[section] 
\newtheorem{lemma}[theorem]{Lemma} 
\newtheorem{proposition}[theorem]{Proposition}

\begin{document}

\baselineskip=15pt

\title[Automorphisms of curves]{Automorphisms of curves
fixing the order two points of the Jacobian}

\author[I. Biswas]{Indranil Biswas}

\address{School of Mathematics, Tata Institute of Fundamental
Research, Homi Bhabha Road, Mumbai 400005, India}

\email{indranil@math.tifr.res.in}

\author{A. J. Parameswaran}

\address{School of Mathematics, Tata Institute of Fundamental
Research, Homi Bhabha Road, Mumbai 400005, India}

\email{param@math.tifr.res.in}

\subjclass[2000]{14H37, 14H40}

\keywords{Curve, automorphism, Jacobian, theta characteristic}

\date{}

\begin{abstract}
Let $X$ be an irreducible smooth projective curve,
of genus at least two, defined
over an algebraically closed field of characteristic
different from two. If $X$ admits a nontrivial
automorphism $\sigma$ that fixes pointwise all the order
two points of $\text{Pic}^0(X)$, then we prove that $X$
is hyperelliptic with $\sigma$ being the unique
hyperelliptic involution. As a corollary, if a nontrivial
automorphisms $\sigma'$ of $X$ fixes pointwise all the
theta characteristics on $X$, then $X$
is hyperelliptic with $\sigma'$ being its hyperelliptic
involution.

\end{abstract}

\maketitle

\section{Introduction}

Let $Y$ be a compact connected Riemann surface of genus
at least two. Assume that there is a nontrivial
holomorphic automorphism
$$
\sigma_0\, :\, Y\, \longrightarrow\, Y
$$
satisfying the condition that for each holomorphic line bundle
$\xi$ over $Y$ with $\xi^{\otimes 2}$ trivializable, the
pull back $\sigma^*_0\xi$ is holomorphically isomorphic to $\xi$.
In \cite{BGP} it was shown that $Y$ must be hyperelliptic and
$\sigma_0$ is the unique hyperelliptic involution (see
\cite[p. 494, Theorem 1.1]{BGP}).

We recall that a theta characteristic on $Y$ is a holomorphic
line bundle $\theta$ such that $\theta^{\otimes 2}$ is
holomorphically isomorphic to the homomorphic cotangent bundle
$K_Y$. The group of order two line bundles on $Y$ acts freely
transitively on the set of all theta characteristics on
$Y$. From this it follows immediately that
if an automorphism of $Y$ fixes pointwise all the
theta characteristics, then it also fixes pointwise all the
order two line bundles on $Y$. Therefore, if $Y$ admits a
nontrivial automorphism $\sigma'_0$ that fixes pointwise
all the theta characteristics on $Y$, then $Y$ is
hyperelliptic and
$\sigma'_0$ is its unique hyperelliptic involution.

The proof of Theorem 1.1 in \cite{BGP} is topological.
Here we investigate the corresponding algebraic geometric
set--up, where the topological proof of Theorem 1.1
in \cite{BGP} is no longer valid.

Let $X$ be an irreducible smooth projective curve defined
over an algebraically closed field $k$. We will assume that
$\text{genus}(X)\, >\, 1$ and $\text{char}(k)\, \not=\, 2$.
We prove the following:

\begin{theorem}\label{thm1}
Let
$$
\sigma\, :\, X\, \longrightarrow\, X
$$
be a nontrivial automorphism that fixes pointwise all the
theta characteristics on $X$. Then $X$ is hyperelliptic
with $\sigma$ being its unique hyperelliptic involution.
\end{theorem}

This theorem is proved by showing that if
$$
\sigma'\, :\, X\, \longrightarrow\, X
$$
is a nontrivial automorphism of $X$ that
fixes pointwise all the
order two points in ${\rm Pic}^0(X)$, then $X$ is
hyperelliptic with
$\sigma'$ being its unique hyperelliptic involution.
(See Lemma \ref{lem1}.)

It should be pointed out that Theorem \ref{thm1} is
not valid if the assumption that the field $k$ is
algebraically closed is removed. There
exists a geometrically irreducible
smooth projective real algebraic curve $Y$ of genus
$g\, \geq\, 2$ which admits a nontrivial involution $\sigma$
that fixes pointwise all the real points
$\xi\, \in\, \text{Pic}^{g-1}(Y)$ with $\xi^{\otimes 2}\,=\,
K_Y$, and $\text{genus}(Y/\langle \sigma\rangle)\,
\not=\,0$. (The details are in \cite{BG}.)

\section{Automorphisms of polarized abelian varieties}

Let $k$ be an algebraically closed field whose characteristic
is different from two. Let $A$ be an abelian variety defined over
$k$ and $L$ an ample line bundle over $A$. For any positive
integer $n$, let
\begin{equation}\label{ker.}
A_n\, \subset\, A
\end{equation}
be the scheme--theoretic kernel of the endomorphism
$A\, \longrightarrow\, A$ defined by $x\, \longmapsto\, nx$.

\begin{proposition}\label{prop1}
Let
$$
\tau\, :\, A\, \longrightarrow\, A
$$
be a nontrivial automorphism such that
$\tau^*L \, =\, L\bigotimes L_0$ for some
$L_0\, \in\, {\rm Pic}^0(A)$, and
the restriction of $\tau$ to the subscheme $A_{n_0}$
(see Eq. \eqref{ker.}) is the identity map for some
$n_0\, \geq\, 2$. Define the two endomorphisms
$$
f_{\pm}\,:=\, {\rm Id}_A\pm \tau\, :\, A\, \longrightarrow\, A\, .
$$
Let $A_+$ (respectively, $A_-$) be the image of $f_+$
(respectively, $f_-$). Then 
\begin{enumerate}
\item $n_0\,=\, 2$.

\item $\tau^2\, =\, \tau\circ\tau$ is the identity automorphism
of $A$.

\item The natural homomorphism
\begin{equation}\label{inc.}
\beta\, :\, A_+ \times A_-\, \longrightarrow\, A
\end{equation}
defined by the inclusions of $A_+$ and $A_-$ in $A$
is an isomorphism.

\item The pull back $\beta^*L$ is of the form $p^*_+L_+
\bigotimes p^*_-L_-$, where $p_+$ (respectively, $p_-$)
is the projection of $A_+ \times A_-$ to
$A_+$ (respectively, $A_-$).
\end{enumerate}
\end{proposition}

\begin{proof}
A proof of statement (1) is given
in \cite[p. 207, Thoerem 5]{Mu}.
See \cite[p. 120, Corollary 1.10]{LB} for a proof under
the assumption that $k$ is the field of complex numbers.

To prove statement (2), we will show that the restriction of
$\tau^2$ to $A_4$ is the identity map. Take any point
$x\,\in\, A_4$. Then $\tau(2x)\,=\,2x$ because
$2x\,\in\, A_2$. Hence
$\tau(x)\,=\, x' + x$ for some $x'\,\in\, A_2$. Thus 
$$
\tau(\tau(x))\,=\,\tau(x' + x)\,=\,\tau(x') +\tau(x) \,=\, x' + 
(x' +x) \,=\, x\, .
$$
Consequently, the restriction of $\tau^2$ to $A_4$ is the identity
map. Now statement (2) follows from statement (1).

To prove statement (3), consider the composition homomorphism
$$
A\,\stackrel{f_+\times f_-}{\longrightarrow}\,
A_+ \times A_- \, \stackrel{\beta}{\longrightarrow}\,  A\, ,
$$
where $\beta$ is the homomorphism in Eq. \eqref{inc.}. It coincides
with the endomorphism of $A$ defined by $x\, \longmapsto\, 2x$.
We also note that $A_2\, \subset\, \text{kernel}(f_+\times f_-)$.
Hence
\begin{equation}\label{kerl.}
\text{kernel}(\beta\circ (f_+\times f_-)) \, \subset\,
\text{kernel}(f_+\times f_-)\, .
\end{equation}

Since $\tau^2\, =\, \text{Id}_A$, the composition $f_+\circ
f_-$ is the zero homomorphism. Hence $\dim (A_+ \times A_{-})\,
\leq\,\dim A$. Now From Eq. \eqref{kerl.} it follows that $\beta$
is an isomorphism.

To prove statement (4), let
$$
\phi_{\beta^*L}\, :\, A_+\times A_-\, \longrightarrow\,
\text{Pic}^0(A_+\times A_-) \, =\, \text{Pic}^0(A_+)\times
\text{Pic}^0(A_-)
$$
be the homomorphism that sends any $k$--rational point $x\, \in\,
A_+\times A_-$
to the line bundle $(t_x^*\beta^*L)\bigotimes \beta^*L^*$, where
$t_x$ is  the
translation map of $A_+\times A_-$ defined by $y\, \longmapsto\,
y+x$; see \cite[p. 131, Corollary 5]{Mu} for a precise definition
of the morphism $\phi_{\beta^*L}$. Let
$$
\tau'\, :=\, \text{Id}_{A_+}\times (-\text{Id}_{A_-})
$$
be the automorphism of $A_+\times A_-$. We note that the
isomorphism $\beta$ in Eq. \eqref{inc.} takes $\tau$ to $\tau'$.

Let
$$
\widehat{\tau}'\, :=\, \text{Id}_{\text{Pic}^0(A_+)}\times 
(-\text{Id}_{\text{Pic}^0(A_-)})
$$
be the automorphism of
$\text{Pic}^0(A_+)\times \text{Pic}^0(A_-)\, =\, \text{Pic}^0(A_+
\times A_-)$. Since $\tau^*L \, =\, L\bigotimes L_0$ for some $L_0
\,\in\, {\rm Pic}^0(A)$, the following diagram is commutative
$$
\begin{matrix}
A_+\times A_- & \stackrel{\phi_{\beta^*L}}{\longrightarrow} &
\text{Pic}^0(A_+)\times \text{Pic}^0(A_-)\\
~\,\Big\downarrow \tau' && ~\,\Big\downarrow \widehat{\tau}'\\
A_+\times A_- & \stackrel{\phi_{\beta^*L}}{\longrightarrow} &
\text{Pic}^0(A_+)\times \text{Pic}^0(A_-)
\end{matrix}
$$
Therefore, the homomorphism $\phi_{\beta^*L}$ takes the
subgroup $A_+$ (respectively, $A_-$) of
$A_+\times A_-$ to the subgroup $\text{Pic}^0(A_+)$
(respectively, $\text{Pic}^0(A_-)$) of $\text{Pic}^0(A_+)\times 
\text{Pic}^0(A_-)$. Now from the injectivity of the homomorphism
$$
\text{NS}(A_+\times A_-)\, \longrightarrow\,
\text{Hom}(A_+\times A_-\, ,\text{Pic}^0(A_+)
\times \text{Pic}^0(A_-))
$$
defined by $\xi\, \longmapsto\,\phi_\xi$ it follows immediately
that the N\'eron--Severi class of $\beta^*L$ coincides with
that of some line bundle of the form 
$p^*_+L_+ \bigotimes p^*_-L_-$ (see \cite[p. 178]{Mu} for the
injectivity of the above homomorphism). Therefore,
statement (4) follows using the fact that
$\text{Pic}^0(A_+)\times \text{Pic}^0(A_-)\, =\, \text{Pic}^0(A_+
\times A_-)$. This completes the proof of the proposition.
\end{proof}

\section{Automorphisms and theta characteristics}

Let $X$ be an irreducible smooth projective
curve, of genus at least two, defined over the field $k$.

\begin{lemma}\label{lem1}
Let
$$
\sigma\, :\, X\, \longrightarrow\, X
$$
be a nontrivial automorphism of $X$ that
fixes pointwise all the
order two points ${\rm Pic}^0(X)_2\, \subset\,
{\rm Pic}^0(X)$. then $X$ is
hyperelliptic with
$\sigma$ being its unique hyperelliptic involution.
\end{lemma}

\begin{proof}
Let $\text{Pic}^d(X)$ denote the moduli space of line
bundles over $X$ of degree $d$. Let $g$ denote the genus of $X$.
On $\text{Pic}^{g-1}(X)$, we have the theta divisor
$\Theta$ given by the
locus of the line bundles admitting nontrivial sections. 
Fix a $k$--rational point $x_0\, \in\, X$. Let $L$ be the
pull back of the line bundle ${\mathcal O}_{\text{Pic}^{g-1}(X)}
(\Theta)$ by the morphism $\text{Pic}^0(X)\, \longrightarrow\, 
\text{Pic}^{g-1}(X)$ that sends any $\zeta$ to $\zeta\bigotimes
{\mathcal O}_X((g-1)x_0)$.

Let $\tau\, :\, \text{Pic}^0(X)\, \longrightarrow\,
\text{Pic}^0(X)$ be the automorphism defined by $\zeta\,
\longmapsto\, \sigma^*\zeta$. This $\tau$ satisfies the
conditions in Proposition \ref{prop1}. Hence $\tau$ is an
involution (see Proposition \ref{prop1}(2)).
This implies that $\sigma$ is an involution.

A hyperelliptic smooth projective curve $Y$ of genus at least
two admits a unique involution $\sigma_Y$ such that
$\text{genus}(Y/\langle\sigma_Y\rangle)\, =\, 0$. Therefore,
to complete the proof of the lemma it suffices to show that
$\text{genus}(X/\langle\sigma\rangle)\, =\, 0$.
We note that the theta divisor $\Theta$ on 
$\text{Pic}^{g-1}(X)$ is irreducible. Indeed, it is the
image of $\text{Sym}^{g-1}(X)$ by the obvious map.
Also, $h^0({\mathcal O}_{\text{Pic}^{g-1}(X)}(\Theta))\, =\,1$
because $\Theta$ defines a principal polarization.

On the other hand, any ample hypersurface of the form
$(A_+\times D_{-})\bigcup (D_+\times A_{-})$ on $A_+\times A_{-}$
is never irreducible unless at least one of $A_+$ and $A_{-}$
is a point; here $D_+$ (respectively, $D_{-}$)
is a hypersurface on $A_+$ (respectively, $A_{-}$).
Therefore, from statement (4) of Proposition \ref{prop1}
and the irreducibility of $\Theta$ we conclude
that either $\dim A_+ \, =\, 0$ or $\dim A_{-} \, =\, 0$.
But $\dim A_{-} \, =\, \text{genus}(X)-
\text{genus}(X/\langle\sigma\rangle)$, and 
$\dim A_{+} \, =\, \text{genus}(X/\langle\sigma\rangle)$.
Since $\text{genus}(X)\, >\, \text{genus}(X/\langle\sigma\rangle)$,
we now conclude that $\text{genus}(X/\langle\sigma\rangle)
\,=\, 0$. This completes the proof of the lemma.
\end{proof}

A line bundle $\theta$ is called 
a \textit{theta characteristic} of $X$ if
$\theta^{\otimes 2}$ is isomorphic to the
canonical line bundle $K_X$ of $X$. The space of theta
characteristics on $X$ is a principal homogeneous space
for ${\rm Pic}^0(X)_2$. Therefore, if an automorphism $\sigma$
of $X$ fixes pointwise all the theta characteristics on $X$,
then $\sigma$ fixes ${\rm Pic}^0(X)_2$ pointwise. Consequently,
the following theorem is deduced from Lemma \ref{lem1}. 

\begin{theorem}\label{thm2}
Let $\sigma\, :\, X\, \longrightarrow\, X$
be a nontrivial automorphism that fixes pointwise all the
theta characteristics on $X$. Then $X$ is hyperelliptic
with $\sigma$ being its unique hyperelliptic involution.
\end{theorem}



\begin{thebibliography}{AAAA}

\bibitem{BG} Biswas, I., Gadgil, S.: Real
theta characteristics and automorphisms of a real curve.
Preprint (2007)

\bibitem{BGP} Biswas, I., Gadgil, S., Sankaran, P.: On
theta characteristics of a compact Riemann surface.
Bull. Sci. Math. \textbf{131}, 493--499 (2007)

\bibitem{LB} Lange, H., Birkenhake, C.: Complex abelian
varieties. Grundlehren der Mathematischen
Wissenschaften, 302. Springer-Verlag, Berlin, 1992

\bibitem{Mu} Mumford, D.: Abelian varieties. Tata Institute
of Fundamental Research Studies in Mathematics, No. 5,
Oxford University Press, London, 1970

\end{thebibliography}
\end{document}